\ifdefined\journalstyle
 \documentclass{baustms}
 \citesort
 \usepackage{latexsym,amsfonts,amsmath,amssymb,url}
 \theoremstyle{cupthm}
 \newtheorem{thm}{Theorem}[section]
 
 \newtheorem{cor}[thm]{Corollary}
 \newtheorem{lemma}[thm]{Lemma}
 \theoremstyle{cupdefn}
 
 \theoremstyle{cuprem}
 
 \numberwithin{equation}{section}

\else
 \documentclass[12pt]{article}
 \usepackage{latexsym,amsfonts,amsmath,amssymb,url}
 \newtheorem{thm}{Theorem}[section]
 \newtheorem{lemma}[thm]{Lemma}
 \newtheorem{cor}[thm]{Corollary}
 \numberwithin{equation}{section}

 \newenvironment{proof}{\begin{trivlist}\item[\hskip\labelsep{\bf Proof.}]}{$\hfill\Box$\end{trivlist}}

 \usepackage[usenames]{color}

\fi

\usepackage{pdfsync}

\newcommand{\satop}[2]{\stackrel{\scriptstyle{#1}}{\scriptstyle{#2}}}

\newcommand{\bsx}{{\boldsymbol{x}}}

\newcommand{\bsX}{{\boldsymbol{X}}}
\newcommand{\rd}{\mathrm{d}}

\newcommand{\bbP}{\mathbb{P}}
\newcommand{\bbR}{\mathbb{R}}

\newcommand{\bbE}{\mathbb{E}}

\newcommand{\bbN}{\mathbb{N}}
\newcommand{\calA}{\mathcal{A}}
\newcommand{\calO}{\mathcal{O}}

\newcommand{\calF}{\mathcal{F}}
\newcommand{\calC}{\mathcal{C}}

\newcommand{\calN}{\mathcal{N}}

\newcommand{\calX}{\mathcal{X}}

\newcommand{\mask}[1]{{}}

\newcommand{\LC}{L\'evy-Ciesielski }
\newcommand{\E}{\mathbb{E}}
\newcommand{\R}{\mathbb{R}}

\definecolor{darkred}{RGB}{139,0,0}
\definecolor{darkgreen}{RGB}{0,100,0}
\definecolor{darkmagenta}{RGB}{139,0,139}
\definecolor{darkpurple}{RGB}{110,0,180}
\definecolor{darkblue}{RGB}{40,0,200}
\definecolor{darkorange}{RGB}{255,140,0}

\definecolor{DARKMAGENTA}{RGB}{139,0,139}  

\begin{document}

\ifdefined\journalstyle

 \title{On the expected uniform error of \\
 Brownian motion approximated by \\the L\'evy-Ciesielski construction}

 \runningtitle{Expected uniform error of L\'evy-Ciesielski construction}

 \author[1]{Bruce Brown}
 \address[1]{School of Mathematics and Statistics, UNSW Sydney 2052, Australia\email{bruce.brown@unsw.edu.au}}

 \author[2]{Michael Griebel}
 \address[2]{Universit\"at Bonn, Institut f\"ur Numerische Simulation,
 Bonn, Germany, and Fraunhofer Institute SCAI, Schloss Birlinghoven,
 Sankt Augustin, Germany\email{griebel@ins.uni-bonn.de}}

 \author[3]{Frances Y.~Kuo}
 \address[3]{School of Mathematics and Statistics, UNSW Sydney 2052, Australia\email{f.kuo@unsw.edu.au}}

 \author[4]{Ian H.~Sloan}
 \address[4]{School of Mathematics and Statistics, UNSW Sydney 2052, Australia\email{i.sloan@unsw.edu.au}}

 \authorheadline{B. Brown, M. Griebel, F. Y. Kuo, and I. H. Sloan}

 \support{The authors acknowledge support of the Australian Research
Council under the project DP210100831. Michael Griebel acknowledges
support from the Sydney Mathematical Research Institute. An earlier
version of this manuscript was posted on arXiv:1706.00915v1.}

\else

 \title{On the expected uniform error of \\
 Brownian motion approximated by \\the L\'evy-Ciesielski construction}

 \author{Bruce Brown, Michael Griebel, Frances Y.~Kuo, and Ian H.~Sloan\footnote{
 School of Mathematics and Statistics, UNSW Sydney 2052, Australia
 ({\tt bruce.brown@unsw.edu.au, f.kuo@unsw.edu.au, i.sloan@unsw.edu.au}) \newline
 Universit\"at Bonn, Institut f\"ur Numerische Simulation,
 Bonn, Germany, and Fraunhofer Institute SCAI, Schloss Birlinghoven,
 Sankt Augustin, Germany ({\tt griebel@ins.uni-bonn.de})}
 \footnote{The authors acknowledge support of the Australian Research
Council under the project DP210100831. Michael Griebel acknowledges
support from the Sydney Mathematical Research Institute.}
 }

 \date{August 2023}

 \maketitle
\fi

\begin{abstract}
It is known that the Brownian bridge or L\'evy-Ciesielski construction of
Brownian paths almost surely converges uniformly to the true Brownian
path. In the present article the focus is on the uniform error. In
particular, we show constructively that at level $N$, at which there are
$d=2^N$ points evaluated on the Brownian path, the uniform error and its
square, and the uniform error of geometric Brownian motion, have upper
bounds of order $\calO(\sqrt{\ln d/d})$, matching the known orders. We
apply the results to an option pricing example.
\end{abstract}

\ifdefined\journalstyle

 \maketitle

\fi

\section{Introduction}

For $t \in [0,1]$, let $B(t) = B(\omega)(t)$ denote the standard Brownian
motion on a probability space $(\Omega,\calF,P)$. That is, for each $t \in
[0,1]$, $B(t)$ is a zero-mean Gaussian random variable, and for each pair
$t,s\in [0,1]$ the covariance is $\bbE[B(t)B(s)] = \min(t,s)$.

In this paper we are concerned with the \LC  (or Brownian bridge)
construction of the Brownian paths. The \LC construction expresses the
Brownian path $B(t)$ in terms of a Faber-Schauder basis $\{ \eta_0,
\eta_{n,i} : n\in\bbN, i=1,\ldots, 2^{n-1}\}$ of continuous functions on
$[0,1]$, where $\eta_0(t) := t$ and
\[
 \eta_{n,i}(t) \,:=\,
 \begin{cases}
 \displaystyle 2^{(n-1)/2} \bigg(t-\frac{2i-2}{2^n}\bigg),
  &\displaystyle t\in \bigg[\frac{2i-2}{2^n}, \frac{2i-1}{2^n}\bigg],\\
 \displaystyle 2^{(n-1)/2} \bigg(\frac{2i}{2^n}-t\bigg),
  &\displaystyle t\in \bigg[\frac{2i-1}{2^n}, \frac{2i}{2^n}\bigg],\\
 0& \mbox{otherwise}.
\end{cases}
\]
For a proof that this is a basis in $\calC[0,1]$, see
\cite[Theorem~2.1(iii)]{Trieb10} or \cite{Trieb10a}. In this construction,
the Brownian path corresponding to the sample point $\omega\in\Omega$ is
given by
\begin{equation}\label{LC}
 B(t) \,:=\,
 X_0(\omega)\, \eta_0(t) +\sum_{n=1}^\infty \sum_{i=1}^{2^{n-1}} X_{n,i}(\omega)\,\eta_{n,i}(t),
\end{equation}
where $X_0$ and all the $X_{n,i}, i=1,\ldots,2^{n-1},n\in \bbN$ are
independent standard normal random variables.
 For $N \in \bbN$ we define the truncated \LC expansion by
\begin{equation}\label{truncatedLC}
 B_{N}(t) \,:=\,
 X_0(\omega)\,\eta_0(t) +
 \sum_{n=1}^N \sum_{i=1}^{2^{n-1}} X_{n,i}(\omega)\,\eta_{n,i}(t).
\end{equation}
Then $B_{N}(t)$ is for each $\omega\in\Omega$ a piecewise-linear function
of $t$ coinciding with $B(t)$ at special values of $t$: we easily see that
$B(0)=B_N(0)=0$, $B(1) =B_N(1) = X_0$, and
$$
 B\bigg(\frac{2\ell-1}{2^N}\bigg)
 \,=\, B_N\bigg(\frac{2\ell-1}{2^N}\bigg), \quad \ell=1,\ldots,2^{N-1},
$$
because the terms in \eqref{LC} with $n> N$ vanish at these points.

The \LC construction has the important property that it converges almost
surely to a continuous Brownian path, see the original works by
\cite{Cie61,Levy} or \cite{Steele}. The precise statement is that, almost
surely,
\[
  \|B-B_N\|_{\infty} \,:=\, \sup_{t\in[0,1]}|B(t) - B_N(t)| \to 0
  \quad\mbox{as}\quad N \to\infty.
\]

The convergence rate for the expected uniform error of the \LC expansion
was obtained in \cite[Theorem~2]{Rit90}: in the language of this paper, we
have
\begin{equation} \label{Ritter}
  \E\, [\|B-B_N\|_{\infty} ] \,\sim\, \sqrt{\frac{\ln d}{2d}},
\end{equation}
where $d := 1+\sum_{n=1}^N 2^{n-1} =2^N$ is the dimension of the
Faber-Schauder basis to level~$N$.

The meaning of the expected value $\bbE$ will be made precise in
Section~\ref{sec:expect}. The asymptotic notation $\alpha(x)\sim \beta(x)$
means that $\lim_{x\to\infty} |\alpha(x)/\beta(x)| \to 1$. Thus
\eqref{Ritter} gives the precise leading term for the expected uniform
error of the \LC expansion. Actually, the article \cite{Rit90} shows also
that the \LC approximation is optimal among all constructions that use
information at $d$ points and Wiener measure. The article \cite{MG02}
contains results for more general problems.

The main result of this paper is Theorem~\ref{thm:main} below which gives
upper bounds of the same order as \eqref{Ritter}, with a slightly worse
constant which is larger by a factor of $2+\sqrt{2} \approx 3.41421$. We
prove this constructively in Section~\ref{sec:B-BN} using a different line
of argument to~\cite{Rit90}, namely extreme value statistics.

\begin{thm}\label{thm:main}
Let $B$ be the \LC expansion of the standard Brownian motion \eqref{LC},
and let $B_N$ be the corresponding truncated expansion
\eqref{truncatedLC}. Then, with $d = 2^N$,
\begin{align*}
 \bbE[\|B-B_N\|_\infty]
 &\,\le\,  (2+\sqrt{2}) \sqrt{\frac{\ln d}{2d}}
\left(1+\calO\left(\frac{1}{\sqrt{\ln d}}\right)\right),
 \\
 \sqrt{\bbE[\|B-B_N\|_\infty^2]}
 &\,\le\,  (2+\sqrt{2}) \sqrt{\frac{\ln d}{2d}}
\left(1+\calO\left(\frac{1}{\sqrt{\ln d}}\right)\right).
\end{align*}
\end{thm}

Geometric Brownian motion is the solution $S(t) = S(\omega)(t)$ at time
$t$ of the stochastic differential equation
\begin{equation} \label{eq:sde}
 \rd S(t) \,=\, S(t) \left(r \,\rd t + \sigma \,\rd B(t)\right),\quad t\in[0,1]
\end{equation}
for given initial data $S(0)$, where $r>0$ is the drift, $\sigma>0$ is the
volatility, and $B(t)$ is the standard Brownian motion. The solution to
\eqref{eq:sde} is given explicitly by
\begin{equation}\label{eq:St}
  S(t) \,=\, S(0)
  \exp\left( \left(r-\tfrac{\sigma^2}{2} \right) t + \sigma B(t)\right).
\end{equation}
Let $S_N$ be the approximation defined by
\begin{equation}\label{eq:SN}
 S_N(t) \,:=\, S(0)
  \exp\left( \left(r-\tfrac{\sigma^2}{2} \right) t + \sigma B_N(t)\right),
\end{equation}
where $B_N$ is the truncated \LC approximation of $B$ given by
\eqref{truncatedLC}. Then we prove in Section~\ref{sec:S-SN} the following
corollary to Theorem~\ref{thm:main}.

\begin{cor} \label{cor:main}
Let $S$ be the geometric Brownian motion \eqref{eq:St}, and let $S_N$ be
the truncated approximation \eqref{eq:SN}. Then, with $d=2^N$,
$$
\E\, [\|S\,-\,S_N \|_\infty]
\,=\, \calO\left(\sqrt{\frac{\ln d}{d}}\right),
$$
where the implied constant depends only on $r$ and $\sigma$.
\end{cor}

Section~\ref{sec:option} gives an application to the problem of pricing an
arithmetic Asian option.

\section{The expected value as an integral over a sequence space}
\label{sec:expect}

In this section we show that the expected value in Theorem~\ref{thm:main}
can be expressed as an integral over a sequence space. We remark that we
will sometimes find it convenient to use interchangeably the language of
measure and integration or alternatively that of probability and
expectation.

Recall that the \LC expansion \eqref{LC} expresses the Brownian path
$B(t)$ in terms of an infinite sequence $\bsX(\omega)=(X_0,
(X_{n,i})_{n\in\bbN, i=1,\ldots, 2^{n-1}})$ of independent standard normal
random variables. In the following we will denote a particular realization
of this sequence $\bsX$ by
\[
  \bsx = \big(x_0,(x_{n,i})_{n\in\bbN, i=1,\ldots,2^{n-1}}\big)
  = (x_0,x_1,x_2,x_3,\ldots) \in\bbR^\infty,
\]
where 
we will switch freely between the double-index
labeling $(x_0,x_{1,1}, x_{2,1}, x_{2,2},\ldots)$ and a single-index
labeling $(x_0,x_1,x_2,x_3,\ldots)$ as appropriate, with the indexing
convention that $x_{n,i}$ becomes $x_{2^{n-1}-1+i}$ for $n\ge 1$ and $1\le
i\le 2^{n-1}$.

It is clear from \eqref{LC} that, for $t\in[0,1]$ and a fixed
$\omega\in\Omega$,
\begin{align} \label{eq:mot}
  \left| B_{N}(t)\right|
  &\,\le\, \left| X_0\right| \,+\,
  \sum_{n=1}^N \left(\max_{1\le i\le 2^{n-1}} \left| X_{n,i}\right|\right)
  \sum_{i=1}^{2^{n-1}}\eta_{n,i}(t)\nonumber\\
  &\,\le\, \left| X_0\right| \,+\,
  \sum_{n=1}^N \max_{1\le i\le 2^{n-1}} \left| X_{n,i}\right| 2^{-(n+1)/2},
\end{align}
where in the last step we used the fact that for a given $n\geq 1$ the
disjoint nature of the Faber-Schauder functions ensures that at most one
value of $i$ contributes to the sum over $i$, and also that the
$\eta_{n,i}$ for $i=1,\ldots,2^{n-1}$ have the same maximum value
$2^{-(n+1)/2}$.

Motivated by the bound \eqref{eq:mot}, and following \cite{GKS_ANOVA}, we
define a norm of the sequence $\bsx = \big(x_0,(x_{n,i})_{n\in\bbN,
i=1,\ldots,2^{n-1}}\big)$ by
\[
 \|\bsx\|_\calX \,:=\, |x_0|+ \sum_{n=1}^\infty \max_{1\leq i \leq 2^{n-1}} |x_{n,i}|\, 2^{-(n+1)/2},
\]
and we define a corresponding normed space by $\calX := \{
\bsx\in\bbR^\infty : \|\bsx\|_\calX < \infty\}$. It is easily seen that
$\calX$ is a Banach space.

Each choice of $\bsx\in\calX$ corresponds to a particular
$\omega\in\Omega$ (but not \emph{vice versa}, since there are sample
points $\omega\in\Omega$ corresponding to sequences $\bsx$ for which the
norm $\|\bsx\|_\calX$ is not finite).  Hence to each $\bsx\in\calX$ there
corresponds a particular Brownian path via \eqref{LC}, or expressed in
terms of $\bsx$,

\begin{equation}\label{eq:hdef}
 B(\bsx)(t) \,=\,x_0\, \eta_0(t) +\sum_{n=1}^\infty \sum_{i=1}^{2^{n-1}} x_{n,i}\,
 \eta_{n,i}(t),\quad t\in[0,1].
\end{equation}
That the resulting path is continuous on $[0,1]$ follows from the
observation that the path is the pointwise limit of the truncated series
\begin{equation}\label{eq:hdefN}
 B_N(\bsx)(t) \,=\,x_0\, \eta_0(t) +\sum_{n=1}^N \sum_{i=1}^{2^{n-1}} x_{n,i}\,
  \eta_{n,i}(t),\quad t\in[0,1],
\end{equation}
which is uniformly convergent since
\[
  \|B_N\|_\infty \,\le\,
  |x_0|+ \sum_{n=1}^\infty
  \max_{1\leq i \leq 2^{n-1}} |x_{n,i}|\, 2^{-(n+1)/2}
  \,=\, \|\bsx\|_\calX \,<\, \infty
  \quad\mbox{for}\quad \bsx\in\calX,
\]
so that \eqref{eq:hdef} does indeed define a continuous function for
$\bsx\in\calX$.

We define $\calA_{\bbR^\infty}$ to be the $\sigma$-algebra generated by
products of Borel sets of $\bbR$, see \cite[p.\ 372]{Bog}. On the Banach
space $\calX$, we now define a product Gaussian measure $\rho(\rd\bsx) :=
\otimes_{j=0}^\infty \phi(x_j) \,\rd x_j$, see \cite[p.\ 392 and Example
2.35]{Bog}, where $\phi$ is the standard normal probability density
$\phi(x) := \exp(- x^2/2)/\sqrt{2\pi}$.

We next show that the space $\calX$ has full Gaussian measure, i.e. that
$$
  \bbP\bigg(|X_0|+ \sum_{n=1}^\infty \max_{1\leq i\leq 2^{n-1}} |X_{n,i}|\,
2^{-(n+1)/2}<\infty\bigg) \,=\,1.
$$
This fact is the basis of the classical proof that the \LC construction
almost surely converges uniformly to the Brownian path. For a brief
explanation, we define
\[
  H_n(\omega) \,:=\,
  \begin{cases}
  \left|X_0(\omega)\right| & \mbox{for } n=0,\\
  \max_{1\leq i\leq 2^{n-1}}\left|X_{n,i}(\omega)\right| 2^{-(n+1)/2} &
  \mbox{for } n \geq 1.
 \end{cases}
\]
As a consequence of the Borel-Cantelli lemma, one can construct a sequence
$(\beta_n)_{n \geq 1}$ of positive numbers such that
$$
 \sum_{n=1}^\infty \beta_n <\infty, \quad\mbox{and}\quad
 \bbP\left(H_n(\cdot) > \beta_n \mbox{ infinitely often}\right) =0.
$$
We now define $\tilde \Omega$ to be the subset of $\Omega$ consisting of
the sample points $\omega$ for which $H_n(\omega) > \beta_n$ for only
finitely many values of $n$. Then $\tilde \Omega$ is of full Gaussian
measure, and for each $\omega \in \tilde \Omega $ there exists $N(\omega)
\in \bbN$ such that $H_n(\omega) \leq \beta_n$ for $n> N(\omega)$, giving
$$
\sum_{n=1}^{\infty} H_n(\omega)
 \,\le\, \sum_{n=1}^{N(\omega)} H_n(\omega) +
 \sum_{n=N(\omega)+1}^\infty \beta_n \,<\, \infty \quad \mbox{for}\;
 \omega \in \tilde \Omega.
$$
Thus $\bbP(\sum_{n=0}^\infty H_n <\infty)=1$ as claimed, and $\calX$ is of
full Gaussian measure.

We now study integration on the measure space
$(\calX,\calA_{\bbR^\infty},\rho)$, and we denote the integral, or the
expected value,  of a measurable function $f$ by $\E[f] := \int_{\calX}
f(\bsx) \,\rho(\rd\bsx)$.

\section{Expected uniform error of standard Brownian motion} \label{sec:B-BN}

We devote this section to proving Theorem~\ref{thm:main}. We have from
\eqref{LC} and \eqref{truncatedLC}
\begin{align*}
 |B(t) - B_N(t)| \,=\,
 \Bigg|\sum_{n=N+1}^\infty \sum_{i=1}^{2^{n-1}} X_{n,i}\, \eta_{n,i}(t) \Bigg|
 \,\le\, \sum_{n=N+1}^\infty \left(\max_{1\le i\le 2^{n-1}} |X_{n,i}|\right)
 \sum_{i=1}^{2^{n-1}} \eta_{n,i}(t).
\end{align*}
Using the same disjoint support argument as in \eqref{eq:mot}, we conclude
that
\begin{align*}
  \|B - B_N\|_\infty
  \,\le\, \sum_{n=N+1}^\infty \max_{1\le i\le 2^{n-1}} |X_{n,i}|\, 2^{-(n+1)/2}
  \,=\, \sum_{\ell = 2^{N}, 2^{N+1},2^{N+2},\ldots} \frac{M_\ell}{2\sqrt{\ell}}\,,
\end{align*}
where we introduced new random variables
\begin{equation}\label{Ml-def}
  M_\ell\,:=\, M_{2^{n-1}} \,:=\, \max_{1\le i \le 2^{n-1}}|X_{n,i}|
\qquad\mbox{for $\ell = 2^{n-1}$  and $n\ge 1$}\,.
\end{equation}
Thus
\begin{align} \label{eq:B2M}
  \bbE\,[ \|B - B_N\|_\infty ]
  \,\le\,
  \sum_{\ell = 2^{N}, 2^{N+1},2^{N+2},\ldots} \frac{\bbE \left[M_\ell \right]}{2\sqrt{\ell}},
\end{align}
and since $M_\ell$ and $M_{\ell'}$ are independent random variables for
$\ell \ne \ell'$,
\begin{align} \label{eq:B2Msq}
  \bbE\,[ \|B - B_N\|_\infty^2 ]
  &\,\le\,
  \sum_{\ell,\ell' = 2^{N}, 2^{N+1},2^{N+2},\ldots}
  \frac{\bbE \left[M_\ell\, M_{\ell'}\right]}{2\sqrt{\ell}\cdot 2\sqrt{\ell'}} \nonumber\\
  &\,=\,
  \sum_{\ell = 2^{N}, 2^{N+1},2^{N+2},\ldots} \frac{\bbE\,\big[M_\ell^2 \big]}{2\sqrt{\ell}\cdot 2\sqrt{\ell}}
  + \sum_{\satop{\ell,\ell' = 2^{N}, 2^{N+1},2^{N+2},\ldots}{\ell\ne\ell'}}
  \frac{\bbE \left[M_\ell\right]\, \bbE\left[M_{\ell'}\right]}{2\sqrt{\ell}\cdot 2\sqrt{\ell'}}.
\end{align}

Now we are in the territory of extreme value statistics.  It is known that
the distribution function of the maximum of the absolute value of $\ell$
independent and identical Gaussian random variables converges (after
appropriate centering and scaling, as below) to the Gumbel distribution. A
first step is to obtain an explicit expression for the distribution
function of $M_\ell$.  Because $X_{n,1}, X_{n,2}, \ldots,X_{n,\ell}$ are
$\calN(0,1)$ random variables, for $x\in\bbR^+$ and $i=1,\ldots,\ell$, we
have
\[
  \bbP(X_{n,i}\le x)
 \,=\,\int_{-\infty}^x \phi(t)\,\rd t\,=:\, \Phi(x),
\]
where $\phi$ is the standard normal density. Similarly,
\[
  \bbP(|X_{n,i}|\le x) \,=\, \int_{-x}^x \phi(t)\,\rd t
  \,=\, \Phi(x) - \Phi(-x)
  \,=\, 2\Phi(x) - 1.
\]
Therefore (since $X_{n,1}, X_{n,2}, \ldots,X_{n,\ell}$ are independent
random variables) we have
\begin{align*}
  \bbP(M_\ell \le x)
  &\,=\, \bbP(|X_{n,1}|\le x \mbox{ and } |X_{n,2}|\le x \mbox{ and }
  \cdots \mbox{and }|X_{n,\ell}|\le x)
  \,=\, (2\Phi(x) - 1)^\ell .
\end{align*}
Thus the distribution function of $M_\ell$ is
\begin{equation}\label{eq:distM}
\Psi_\ell(x) \,:=\, (2\Phi(x) - 1)^\ell,\quad x\in \mathbb{R^+}.
\end{equation}

We now define a new random variable $Y_\ell$ for $\ell\ge 1$, which is a
recentered and rescaled version of $M_\ell$:
\begin{equation}\label{Yl}
 Y_\ell \,:=\, \frac{M_\ell-a_\ell}{b_\ell},\quad
 \mbox{or equivalently,}\quad
  M_\ell \,=\, a_\ell + b_\ell Y_\ell, \quad a_\ell >0,\, b_\ell >0.
\end{equation}
It is known (see below) to be appropriate to take $a_\ell$ and $b_\ell$ to
satisfy
\begin{equation}\label{al}
 a_\ell =\sqrt{2\ln \ell}+o(1), \quad\mbox{and}\quad
 b_\ell=\frac{1}{a_\ell}.
\end{equation}
More precisely, for later convenience we will define $a_\ell$ to be the
unique solution of
\begin{equation}\label{eq:aldef}
 \frac{1}{\ell} = \sqrt{\frac{2}{\pi}}\frac{e^{-a_\ell^2/2}}{a_\ell} =: g(a_\ell).
\end{equation}
We now show that \eqref{eq:aldef} implies \eqref{al}.

\begin{lemma} \label{lem:al}
Equation \eqref{eq:aldef} for $\ell\ge 1$ has a unique positive solution
of the form $a_\ell=\sqrt{2\ln{\ell}}+o(1)$. Moreover, for $\ell\ge 3$ we
have $a_\ell \in (1,\sqrt{2\ln \ell})$.
\end{lemma}

\begin{proof}
The fact that any solution of \eqref{eq:aldef} is positive is immediate.
Now observe that $g$ in \eqref{eq:aldef} is monotonically decreasing from
$+\infty$ to $0$ on $\bbR^+$. It follows immediately that there is a
unique solution $a_\ell \in (0,\infty)$ for \eqref{eq:aldef}. Moreover, we
have
\[
 a_\ell > 1
 \quad\Leftrightarrow\quad
 \frac{1}{\ell}
 < \sqrt {\frac 2 \pi }\frac {e^{-1^2/2}} {1}
 = \sqrt {\frac 2 {\pi e} } =0.484\ldots,
\]
which holds if and only if $\ell\ge 3$. Now observe that \eqref{eq:aldef}
is equivalent to
\begin{equation}\label{eq:aelleq}
 a_\ell = \sqrt{ 2\left(\ln \ell - \ln \left(\sqrt {\frac  \pi 2} a_\ell\right)\right)}.
\end{equation}
For $\ell \ge 3$ we have $a_\ell > 1$ and hence $\ln(\sqrt{\pi/2}\,a_\ell)>\ln(\sqrt{\pi/2})>0$,
so from \eqref{eq:aelleq} we have $a_\ell < \sqrt{ 2 \ln \ell }$. In turn it follows that
\[
  a_\ell > \sqrt{2 \ln \ell - 2 \ln\left(\sqrt {\frac \pi 2 }\sqrt{2 \ln \ell}\right)}.
\]
Thus for $\ell\ge 3$ we have $1\le a_\ell=\sqrt{2\ln \ell} +o(1) \le
\sqrt{2\ln{\ell}}$, completing the proof.
\end{proof}

It is well known that the distribution function of $Y_\ell$ converges in
distribution to a random variable with the Gumbel distribution
$\exp(-e^{-y})$. For later convenience we state this as a lemma and give a
short proof.

\begin{lemma} \label{lem:Gumbel}
The random variable $Y_ \ell$ defined in \eqref{Yl}, with $a_\ell$ defined
by \eqref{eq:aldef} and $b_\ell = 1/a_\ell$, converges in distribution to
a random variable $Y$ with Gumbel distribution function $\bbP(Y\le y) =
\exp(-e^{-y})$.
\end{lemma}

\begin{proof}
The proof is based on the asymptotic version of Mill's ratio \cite{Small},
\[
1-\Phi(x)\,\sim\, \frac{\phi(x)}{x},\quad x\to +\infty,
\]
where, as in the Introduction, $\sim$ means that the quotient of the two
sides converges to~$1$. From this it follows that for $y\in\R$
\begin{align*}
\bbP(Y_\ell\le y)
&\,=\,\bbP(M_\ell\le a_\ell+b_\ell y)\,=\,(2\Phi(a_\ell+b_\ell y)-1)^\ell
\,=\,(1-2[1-\Phi(a_\ell+b_\ell y)])^\ell\\
&\,\sim\,\bigg(1-\sqrt{\frac{2}{\pi}}
\frac{\exp(-\tfrac{1}{2}(a_\ell  +  b_\ell y)^2)}{a_\ell+b_\ell y}\bigg)^\ell
\,\sim\,\bigg(1-\sqrt{\frac{2}{\pi}}
\frac{\exp(-\tfrac{1}{2}a_\ell^2 - a_\ell b_\ell y)}{a_\ell}\bigg)^\ell \\
&\,=\,\big(1- \exp(-y)/\ell\big)^\ell
\,\sim\,\exp(-e^{-y}) \quad \mbox{as} \;\ell\to\infty,
\end{align*}
where in the second step we dropped a higher order term, and in the second
last step we used \eqref{eq:aldef} and $a_\ell b_\ell = 1$, thus proving
the lemma.
\end{proof}

A deeper result, which we need, is that $Y_\ell$ converges in expectation
to the limit~$Y$. This is proved in the following lemma.

\begin{lemma} \label{lem:Yl}
The random variable $Y_\ell$ defined in \eqref{Yl}, with $a_\ell$ defined
by \eqref{eq:aldef} and $b_\ell = 1/a_\ell$, converges in expectation to a
random variable~$Y$ with Gumbel distribution $\exp(-e^{-y})$, thus
\[
\lim_{\ell\to\infty}\mathbb{E}[Y_\ell] \,=\, \bbE[Y]
\,=\, \int_{-\infty}^\infty y \exp(-y-e^{-y})\rd y \,=\, \gamma,
\]
where $\gamma$ is Euler's constant.
\end{lemma}

\begin{proof}
For a sequence of real-valued random variables $Y_1,Y_2,\ldots$ converging
in distribution to a random variable $Y$, it is well known that a
sufficient condition for convergence in expectation is uniform
integrability of the $Y_\ell$. In turn a sufficient condition for uniform
integrability is that for sufficiently large $\ell$
\begin{align*}
&\bbP(Y_\ell \ge y) \le Q(y) \mbox{ for } y > 0, \qquad\mbox{and}\qquad
\bbP(Y_\ell \le y) \le R(y) \mbox{ for } y < 0,
\end{align*}
where $Q(y)$ is integrable on $\R^+$ and $R(y)$ is integrable on $\R^-$.

First assume $y>0$. We have from \eqref{eq:distM} that
\begin{align*}
 \bbP(Y_\ell\ge y)
 &\,=\, \bbP(M_\ell\ge a_\ell + b_\ell y)
 \,=\, 1- \bbP(M_\ell\le a_\ell + b_\ell y) \\
 &\,=\, 1- (2\Phi(a_\ell + b_\ell y)-1)^\ell
 \,=\, 1- (1 - 2[1-\Phi(a_\ell + b_\ell y)])^\ell \\
 &\,\le 1-\left(1-2\frac{\phi(a_\ell + b_\ell y)}{a_\ell + b_\ell y}\right)^\ell
 \,=\,1-\left(1-\sqrt{\frac{2}{\pi}}
 \frac{\exp(-\tfrac{1}{2}(a_\ell + b_\ell y)^2)}{a_\ell + b_\ell y}\right)^\ell\\
 &\,\le\, 1-\left(1-\sqrt{\frac{2}{\pi}}\frac{\exp\left(-\tfrac{1}{2}a_\ell^2 - a_\ell b_\ell y\right)}
 {a_\ell}\right)^\ell,
\end{align*}
where we used the upper bound form of Mills' ratio \cite{Small},
\begin{equation} \label{eq:MR-upper}
1-\Phi(x)< \frac{\phi(x)}{x},\quad x\in\mathbb{R}^+,
\end{equation}
and dropped harmless terms in both the denominator and the exponent in the
numerator. Using now \eqref{eq:aldef} and also $a_\ell b_\ell= 1$ we have
\begin{align} \label{eq:defQ}
 \bbP(Y_\ell\ge y)
 &\,\le\, 1-\left(1-\frac{\exp(-y)}{\ell}\right)^\ell
 \,\le\, \exp(-y) \,=:\, Q(y)\,,
\end{align}
where we used the fact that the function $(1-c/x)^x$ is increasing on
$[1,\infty)$ for $c\in [0,1]$, and hence takes its minimum at $x=1$. It
follows that
\[
 \int_0^\infty \bbP(Y_\ell\ge y)\,\rd y \,\le\, \int_0^\infty
 \exp(-y)\,\rd y \,=\, 1.
\]

Now we consider $y<0$. Note first that $M_\ell=a_\ell+b_\ell Y_\ell$ takes
only non-negative values, thus we may restrict $y$ to $y\ge
-a_\ell/b_\ell$. We have
\begin{align*}
 \bbP(Y_\ell\le y)
 &\,=\,\bbP(M_\ell\le a_\ell + b_\ell y)
 \,=\, (2\Phi(a_\ell + b_\ell y)-1)^\ell\,.
\end{align*}
Now for $t>0$ the standard normal distribution $\Phi$ has negative second
derivative,
$$
\Phi''(t)=\phi'(t) <0\quad\mbox{for } t>0,
$$
and first derivative $\Phi'(t)=\phi(t)$, from which it follows that
\begin{align*}
\Phi(a_\ell + b_\ell y)\,&\le\, \Phi(a_\ell) + b_\ell y\,\phi(a_\ell)
\quad \mbox{for } y\ge -a_\ell/b_\ell.
\end{align*}
Thus on using $b_\ell=1/a_\ell$, we obtain
\begin{align*}
\bbP(Y_\ell\le y)
&\,\le\,\bigg(2\Phi(a_\ell)+2a_\ell^{-1}y\phi(a_\ell)-1\bigg)^\ell\\
&\,\le\,\bigg(1-2a_\ell^{-1}\phi(a_\ell)(1-y-a_\ell^{-2})\bigg)^\ell
\,=\,\bigg(1-\frac{1}{\ell}(1-y-a_\ell^{-2})\bigg)^\ell,
\end{align*}
where in the second last step we used the lower bound form of Mills'
ratio, see \cite[p.~44]{Small}
\begin{equation} \label{eq:MR-lower}
1-\Phi(t)\ge \frac{\phi(t)}{t}(1-t^{-2})\quad\mbox{for } t>0,
\end{equation}
and in the last step we used \eqref{eq:aldef}.  If we now take $\ell\ge L$
then we have
\begin{align} \label{eq:defR}
 \bbP(Y_\ell\le y)
 &\,\le\,\bigg(1-\frac{1}{\ell}(1-y-a_L^{-2})\bigg)^\ell \nonumber\\
 &\,\le\, \exp\left(-(1-y-a_L^{-2})\right)
 \,=\,\exp\left(-(1-a_L^{-2})\right) \exp(y)
 \,=:\,R(y),
\end{align}
since the convergence in the last limit is monotone increasing.  The
function $R(y)$ so defined is integrable on $\R^-$, completing the proof
that $Y_\ell$ converges in expectation.

It then follows from Lemma~\ref{lem:Gumbel} that the limit of
$\bbE[Y_\ell]$ is precisely $\bbE[Y] = \gamma$.
\end{proof}

Since Lemma~\ref{lem:Yl} establishes the convergence of $\bbE[Y_\ell]$ as
$\ell \to \infty$, it can be inferred that there exists a positive
constant $c$ such that
\begin{align} \label{eq:def-c}
  \bbE[Y_\ell]\le c
  \quad\mbox{and hence}\quad
  \bbE[M_\ell] \,\le\, a_\ell + b_\ell\,c\, \le  a_\ell + c,
\end{align}
where we used $b_\ell = a_\ell^{-1} \le 1$ for $\ell\ge 3$. We then
conclude from \eqref{eq:B2M} that
\begin{align}\label{EBmBN}
  \bbE\,[ \|B - B_N\|_\infty ]
  &\,\le\, \sum_{\ell = 2^{N}, 2^{N+1},2^{N+2},\ldots} \frac{a_\ell+c}{2\sqrt{\ell}}.
\end{align}

It only remains to estimate the sum in \eqref{EBmBN}. Using
Lemma~\ref{lem:al} with $N\ge 2$ (and hence $\ell\ge 3$), we have $a_\ell
< \sqrt{2\ln \ell}$, and on setting $\ell= 2^{N+j}$,
\begin{align*}
 \sum_{\ell=2^N, 2^{N+1}, 2^{N+2},\cdots}\frac{a_\ell}{\sqrt{\ell}}
 &\,\le\,
 \sum_{j=0}^\infty \frac {\sqrt{2(\ln 2) (N+j)}}{\sqrt{2^{N+j}}}
 \,=\,  2^{-(N-1)/2} \sqrt{\ln 2}\,	\sum_{j=0}^\infty \frac {\sqrt{ N+j}}{2^{j/2}} .\nonumber \\
 &\,\le\, 2^{-(N-1)/2} \sqrt{\ln 2}\,\sqrt{N}\,(2+\sqrt 2 )\,(1 + \calO(N^{-1/2})), \nonumber
\end{align*}
where in the final step we used $\sqrt{N+j} \leq \sqrt{N}+\sqrt{j}$ and
$\sum_{j=0}^\infty 1/2^{j/2}=2+\sqrt{2}$, while noting that
$\sum_{j=0}^\infty \sqrt{j}/2^{j/2}$ is finite and independent of $N$.
Moreover, by a similar argument we conclude that $\sum_{\ell=2^N, 2^{N+1},
2^{N+2},\cdots} c/\sqrt{\ell} = \calO(2^{-N/2})$, thus altogether we
obtain from \eqref{EBmBN}
\begin{align*}
 &\bbE\,[ \|B - B_N\|_\infty ]
 \,\le\, \frac{1}{2} \cdot 2^{-(N-1)/2}\sqrt{\ln2}\,\sqrt{N}\,(2+\sqrt{2})\,(1+\calO(N^{-1/2})) \\
 &\,=\, (2+\sqrt 2) \frac {\sqrt{N\ln{2}}}{\sqrt {2\cdot 2^{N}}}(1+\calO(N^{-1/2}))
 \,=\, (2+\sqrt 2) \frac {\sqrt{\ln{d}}}{\sqrt {2 d}}\left
(1+\calO\left(\frac{1}{\sqrt{\ln{d}}}\right)\right),
\end{align*}
which proves the first bound in Theorem~\ref{thm:main}.

To prove the second bound in Theorem~\ref{thm:main}, we need first to
bound $\bbE\,\big[M_\ell^2\big]$. With $Y_\ell$ and $Y$ defined as above,
for $v>0$ we have
\begin{align} \label{eq:Yl2}
 \bbP(Y_\ell^2 \ge v)
 &\,=\, \bbP(Y_\ell \ge \sqrt{v}) + \bbP(Y_\ell \le - \sqrt{v}) \\
 &\to\, \bbP(Y \ge \sqrt{v}) + \bbP(Y \le - \sqrt{v})
 \,=\, \bbP(Y^2 \ge v)
 \quad\mbox{as}\quad \ell\to\infty, \nonumber
\end{align}
while for $v<0$ we have $\bbP(Y_\ell^2 \ge v) = \bbP(Y^2 \ge v)= 1$. Thus
by Lemma~\ref{lem:Gumbel}, $Y_\ell^2$ converges in distribution to $Y^2$,
where $Y$ is the Gumbel distribution. To prove convergence in expectation,
we use \eqref{eq:Yl2} with \eqref{eq:defQ} and \eqref{eq:defR} to give,
for $\ell\ge L$ and $v>0$,
\[
  \bbP(Y_\ell^2 \ge v) \,\le\, \exp(-\sqrt{v}) + \exp\left(- (1-a_L^{-2})\right) \exp(-\sqrt{v}),
\]
which is integrable on $\bbR^+$, proving $\bbE\,\big[Y_\ell^2\big] \to
\bbE\,\big[Y^2\big] < \infty$. In turn it follows that there exists $c'>0$
such that $\bbE\,\big[Y_\ell^2\big] \le c'$, and together with
\eqref{eq:def-c} we have
\[
  \bbE\,\big[M_\ell^2\big]
  \,\le\, a_\ell^2 + 2\,a_\ell\, b_\ell\, c + b_\ell^2\, c'
  \,\le\, a_\ell^2 + 2\,a_\ell\,c + c'
  \,\le\, (a_\ell + c'')^2,
\]
where we used $b_\ell = a_{\ell}^{-1} \le 1$ for $\ell\ge 3$ and
introduced $c'' := \max(c,\sqrt{c'})$. We now have from \eqref{eq:B2Msq}
that
\begin{align*}
  \bbE\,[ \|B - B_N\|_\infty^2 ]
  &\,\le\,
  \sum_{\ell = 2^{N}, 2^{N+1},2^{N+2},\ldots} \frac{(a_\ell + c'')^2}{2\sqrt{\ell}\cdot 2\sqrt{\ell}}
  + \sum_{\satop{\ell,\ell' = 2^{N}, 2^{N+1},2^{N+2},\ldots}{\ell\ne\ell'}}
  \frac{(a_\ell + c'')(a_{\ell'} + c'')}{2\sqrt{\ell}\cdot 2\sqrt{\ell'}} \\
  &\,=\, \left(\sum_{\ell = 2^{N}, 2^{N+1},2^{N+2},\ldots} \frac{a_\ell + c''}{2\sqrt{\ell}}\right)^2,
\end{align*}
which is the square of the right-hand side of \eqref{EBmBN}, with $c$
replaced by $c''$. The second bound in Theorem~\ref{thm:main} then
follows.

\section{Expected uniform error of geometric Brownian motion}
\label{sec:S-SN}

We are now in the position to give a proof of Corollary~\ref{cor:main}.
{}From \eqref{eq:St} and \eqref{eq:SN} it follows that $S(t) - S_N(t)
  \,=\, S(0)\, e^{(r-\sigma^2/2)t}\,
  \big(\exp(\sigma B(t)) -\exp(\sigma B_N(t))\big)$,
and thus
\begin{align} \label{eq:extrastep}
  \|S - S_N\|_\infty
  &\,\le\, S(0)\,e^{|r-\sigma^2/2|}\,
  \big\|\exp(\sigma B) -\exp(\sigma B_N)\big\|_\infty.
\end{align}
In turn it follows that
\begin{align} \label{eq:same}
  \|S - S_N\|_\infty
&\,\le\,S(0)\,e^{|r-\sigma^2/2|}\,
  \big\|\exp(\sigma B_N) \big(\exp(\sigma (B-B_N)) - 1\big)\big\|_\infty
  \nonumber\\
&\,\le\,S(0)\,e^{|r-\sigma^2/2|}\,
  \exp\big(\sigma \|B_N\|_\infty\big) \big(\exp(\sigma \|B-B_N\|_\infty) - 1\big)  .
\end{align}

Using $|\exp(x)-1| \le |x|\,\exp(|x|)$ for $x\in\bbR$ and $\|B_N\|_\infty
\le \|B\|_\infty$ for $N\in\bbN$, we have
\begin{align*}
  \|S - S_N\|_\infty
&\,\le\,S(0)\,e^{|r-\sigma^2/2|}\,
  \exp\big(\sigma \|B_N\|_\infty\big)\,\sigma\, \|B-B_N\|_\infty\,\exp(\sigma \|B-B_N\|_\infty)  \\
&\,\le\,S(0)\,e^{|r-\sigma^2/2|}\,\sigma\,
  \exp\big(3\sigma \|B\|_\infty\big)\, \|B-B_N\|_\infty,
\end{align*}
where we used $\|B-B_N\|_\infty \le \|B\|_\infty + \|B_N\|_\infty \le
2\|B\|_\infty$. By the Cauchy-Schwarz inequality we obtain
\begin{align*}
  \bbE \left[ \|S - S_N\|_\infty \right]
&\,\le\,S(0)\,e^{|r-\sigma^2/2|}\,\sigma\,
  \sqrt{\bbE \left[\exp\big(6\sigma \|B\|_\infty\big)\right]}\,
  \sqrt{\bbE \left[\|B-B_N\|_\infty^2 \right]}.
\end{align*}
It is well known that $\bbE \left[\exp\big(\alpha \|B\|_\infty\big)\right]
< \infty$ for every $\alpha> 0$. Hence the result now follows from the
second bound of Theorem~\ref{thm:main}.

\section{Application to option pricing} \label{sec:option}

Now we consider a continuous version of a path-dependent call option with
strike price $K$ in a Black-Scholes model with risk-free interest rate
$r>0$ and constant volatility $\sigma>0$. Recall that the asset price
$S(t)$ at time $t$ is given explicitly by \eqref{eq:St}. The discounted
payoff for the case of a continuous arithmetic Asian option with terminal
time $T=1$ is therefore
\begin{align}\label{eq:P}
  P \,:=\,&\, e^{-rT}\max\left(\frac 1 T \int_0^T S(t) \,\rd t -K,0\right) \nonumber\\
  \,=\,&\, e^{-r}\max\bigg( S(0)\int_0^1
  \exp\Big(\big(r-\tfrac{\sigma^2}{2}\big)\,t+\sigma B(t) \Big)\,\rd t - K, 0 \bigg).
\end{align}
The pricing problem is then to compute the expected value $\bbE(P)$.

We use the \LC expansion \eqref{eq:hdef} and \eqref{eq:hdefN} for $B(t)$
and $B_N(t)$, and define
\begin{align}\label{eq:PN}
  P_N
  \,:=\, e^{-r}\max\bigg( S(0)\int_0^1
  \exp\Big(\big(r-\tfrac{\sigma^2}{2}\big)\,t+\sigma B_N(t) \Big)\,\rd t - K, 0 \bigg).
\end{align}
We are interested in estimating how fast $\bbE\, [|P - P_N|]$ converges to
$0$ as $N\to\infty$.

\begin{cor}
For $P$ and $P_N$ defined by \eqref{eq:P} and \eqref{eq:PN}, we have
\[
 \bbE\, [|P - P_N|]
  \,=\, \calO\left(\sqrt{\frac{\ln d}{d}}\right),
\]
where $d=2^N$ and the implied constant is independent of $d$.
\end{cor}

\begin{proof}
It can be easily verified that
\begin{equation*}
  |\max(\alpha-K,0) - \max(\beta-K,0)| \,\le\, |\alpha-\beta|.
\end{equation*}
Thus
\begin{align*}
  |P- P_N|
  &\,\le\, \bigg| e^{-r} S(0)\int_0^1 e^{(r-\sigma^2/2)t}
  \Big(\exp(\sigma B(t)) -\exp(\sigma B_N(t))\Big)
  \,\rd t\, \bigg| \\
  &\,\le\, e^{-r} S(0)\,e^{|r-\sigma^2/2|} \int_0^1
  \Big|\exp(\sigma B(t)) -\exp(\sigma B_N(t))\Big|
  \,\rd t \\
  &\,\le\, e^{-r} S(0)\,e^{|r-\sigma^2/2|}
  \|\exp(\sigma B) - \exp(\sigma B_N)\|_\infty,
\end{align*}
where the last upper bound differs from the upper bound
\eqref{eq:extrastep} on $\|S-S_N\|_\infty$ only by a factor of $e^{-r}$.
Hence the result follows from Corollary~\ref{cor:main}.
\end{proof}

%
%
%
%
%
%
%
%
%
\end{document}